\documentclass{article}[12pt]
\usepackage{amsmath,amsfonts,amssymb}
\usepackage[T2A]{fontenc}
\usepackage[cp866]{inputenc}
\usepackage[russian]{babel}
\newcommand{\la}{\lambda}
\begin{document}

\large

$$
$$

\begin{center}
\bf LANDAU DISPERSION RELATIONSHIP \\
IN SELF-CONSISTENT FIELD THEORY
\end{center}

\medskip

\begin{center}
\bf S.A.Stepin
\end{center}

\bigskip
\smallskip

{\bf Abstract} \, The spectral problem is studied associated with
Maxwell-Boltzmann equations describing collisionless plasma.
Formula for instability index is obtained and effective conditions
of two-stream instability are given.

\bigskip
\smallskip

\begin{center}
\bf \S\,1. Maxwell-Boltzmann equations of self-consistent field
\end{center}

\bigskip
\smallskip

The main object of our considerations will be the system of
equations describing collisionless plasma consisting of two types
of the particles --- electrons and ions in the presence of
electromagnetic field. Distribution of the particles in plasma is
characterized by the corresponding densities (distribution
functions) $\,f_e\,$ and $\,f_i\,$ dependent on time $\,t,\,$
space coordinate $\,x\,$ and velocity $\,v\,$ of the particles.
Thus the system under consideration is composed of Maxwell
equations for the components of electromagnetic field and kinetic
(collisionless) Boltzmann equations for distribution functions. In
spatially one-dimensional case magnetic induction is trivial
$\,B(t,x)\equiv 0\,$ and the system of equations in question has
the form
\begin{eqnarray*}
\frac{\partial f_e}{\partial t}&+&v\,\frac{\partial f_e}{\partial
x}\,\,\,-\,\,\,\frac{e}{m}\,E\,\frac{\partial f_e}{\partial
v}\,\,=\,\, 0\,,\\
\frac{\partial f_i}{\partial t}&+&v\,\frac{\partial f_i}{\partial
x}\,\,\,+\,\,\,\frac{Ze}{M}\,E\,\frac{\partial f_i}{\partial
v}\,\,=\,\, 0\,,\\
\frac{\partial E}{\partial t}&+&\,4\pi
e\int_{-\infty}^{\infty}v\big(Zf_i-f_e\big)\,dv\,\,=\,\,0\,,\\
\frac{\partial E}{\partial x}&-&4\pi
e\int_{-\infty}^{\infty}\big(Zf_i-f_e\big)\,dv\,\,=\,\,0\,,
\end{eqnarray*}
where $\,e\,$ and  $\,Ze\,$ stand for the charges, $\,m\,$ and
$\,M\,$ are the masses of electrons and ions respectively, while
$\,E=E(t,x)\,$ denotes the electric field strength. Each of the
kinetic equations admits representation in Liouvillean form and
expresses the fact that the total derivative of distribution
function vanishes along the trajectory for the particle of the
proper type.

Henceforth a solution to the stationary system such that
$$
f_e=f_{0e}(v)\,,\quad f_i=f_{0i}(v)\,,\quad E_0(x)\equiv 0\,
$$
is chosen as an unperturbed one. Following [1] we will assume that
ions distribution function $\,f_{0i}(v)\,$ is fixed and linearize
the system with respect to the stationary solution
$\,f_{0e}(v),\,$ $E_0(x)\equiv 0.\,$ As a result we obtain the
system
\begin{eqnarray}
\frac{\partial f}{\partial t}&=&-\,v\,\frac{\partial f}{\partial x}\,\,+\,\,\frac{e}{m}\,f_0'(v)\,E\,,
\label{form1}\\
\frac{\partial E}{\partial t}&=&4\pi
e\int_{-\infty}^{\infty}v\,f\,dv
\label{form2}
\end{eqnarray}
where $\,f=f(t,x,v)\,$ is the perturbation of electrons stationary
distribution function $\,f_0(v)$ $=f_{0e}(v)\,$ while $\,E(t,x)\,$
stands for electric field strength.

Electromagnetic field induced by the electrons traffic and in turn
affecting the evolution of their density is known to be called
self-consistent. It is assumed that distribution functions of the
particles and components of electromagnetic field decay at
infinity sufficiently rapidly. According to charge conservation
law one has
$$
\langle f\rangle\,\,=\,\iint f(t,x,v)\,dx\,dv\,\,=\,\,0\,.
$$
Note that equation
\begin{equation}
\frac{\partial E}{\partial x}\,\,=\,\,-\,4\pi
e\int_{-\infty}^{\infty}f\,dv
\label{form3}
\end{equation}
is compatible with the system (\ref{form1})-(\ref{form2}) under
consideration in the sense that equation (\ref{form1}) integrated
with respect to velocity variable
$$
\frac{\partial}{\partial t}\int_{-\infty}^{\infty}
f\,dv\,\,+\,\,\frac{\partial}{\partial x}\int_{-\infty}^{\infty}
v\,f\,dv\,\,=\,\,0\,,
$$
by virtue of (\ref{form2}) and (\ref{form3}) becomes an identity
$\,E_{tx}\,=\,E_{xt}\,.$

Equation (\ref{form1}) treated as inhomogeneous one with respect
to $\,f(t)\,$ can be reduced to Duhamel type integral equation
\begin{equation}
f(t)\,\,=\,\,U(t)\,f^{(0)}\,+\,\,\frac{e}{m}\,f_0'(v)\int_0^tU(t-s)\,E(s)\,ds
\label{form4}
\end{equation}
where the group of translations
$\,U(t)\!:\varphi(x,v)\mapsto\varphi(x-vt,v)\,$ specifies
evolution generated by homogeneous equation
$\,\,\partial_tf\,+\,v\,\partial_xf\,=\,0\,\,$ while initial
perturbation of distribution function
$\,f\big|_{t=0}=\,\,f^{(0)}\,\,$ satisfies zero mean condition
$\,\langle f^{(0)}\rangle\,=\,0.$

\medskip

Substitution of $f(t)$ expressed from (\ref{form4}) into equation
(\ref{form2}) implies integro-differential equation for the
strength $\,E(t,x)\,$ subject to condition
$\,E\big|_{t=-\infty}\!=\,0\,$ which agrees with the original
setting of the problem about stability or instability of plasma
oscillations. Provided that initial perturbation of the density
$\,f^{(0)}\,$ is symmetric with respect to velocity
$\,f^{(0)}(x,v)=f^{(0)}(x,-v)\,$ while stationary (unperturbed)
distribution function $\,f_0(v)\,$ satisfies condition
$\,f_0(0)=0\,$ the strength of electric field $\,E=E(t,x)\,$
proves to be a solution to Fredholm type equation
\begin{equation}
E\,\,\,=\,\,\,E^{(0)}\,+\,\,\,K\,E\,
\label{form5}
\end{equation}
where $\,\,\displaystyle{E^{(0)}(t,x)\,\,=\,\,-\,4\pi
e\int_{-\infty}^xdy\int_{-\infty}^{\infty}}f^{(0)}(y-vt,v)\,dv\,\,$
with integral operator $\,K\,$ given by the formula
$$
K:\,\,E(t,x)\,\,\mapsto\,\,\,\frac{4\pi
e^2}{m}\int_{-\infty}^{\infty}f_0(v)\,dv\int_0^ts\,E(t-s,x-vs)\,ds\,.
$$

The question about existence and uniqueness of solution to
equation (\ref{form5}) in different functional spaces is reduced
in [1] to evaluation of corresponding norms for integral operator
$\,K\,$ while solution itself is given by perturbation theory
series known as Neumann expansion
$$
E\,\,\,=\,\,\,E^{(0)}\,\,+\,\,\,K\,E^{(0)}\,\,+\,\,K^2\,E^{(0)}\,\,+\,\,\ldots\,.
$$
Another approach to investigation of solvability problem for
equation (\ref{form5}) takes advantage of Laplace-Fourier
transformation in variables $\,t\,$ and $\,x.\,$ After the passage
to associated representation the solvability condition is
formulated in terms of the corresponding Fredholm resolvent
denominator
$$
\Delta(k,\lambda)\,\,=\,\,1\,\,+\,\,\frac{4\pi
e^2}{m}\int_{-\infty}^{\infty}\frac{f_0(v)} {(\lambda +
ikv)^2}\,\,dv
$$
where $\,\lambda\in\mathbb C\,$ is the spectral parameter while
$\,k\in\mathbb R\,$ is the wave number in $\,x$-axis direction.
Namely condition $\,\,\Delta(k,\lambda)\ne 0\,\,$ implies
invertibility of operator $\,I-K.\,$ In due turn the roots of the
so called Landau dispersion relationship
\begin{equation}
\,\,\Delta(k,\lambda)\,\,=\,\,1\,\,+\,\,\frac{4\pi
e^2}{ikm}\int_{-\infty}^{\infty}\frac{f_0'(v)} {\lambda +
ikv}\,\,dv\,\,=\,\,0\,\, \label{form6}
\end{equation}
give rise to solutions of harmonic wave type depending on $\,t\,$
and $\,x\,$ exponentially $\,\,\exp(\lambda t+ikx)\,\,$ which
induce unstable eigenmode regimes (undamped oscillations) in
plasma provided that $\,{\rm Re}\,\la>0.\,$ In what follows for
notational convenience we will assume that $\,4\pi e^2/m=1.$

\bigskip
\bigskip

\begin{center}
\bf \S\,2. Spectral problem and Schur complement
\end{center}

\bigskip
\medskip

Linear dynamical system (\ref{form1})-(\ref{form2}) under
consideration is associated with infinitesimal operator $\,T\,$
given by the formula
$$
T\,:\,{\displaystyle\begin{pmatrix}
f(x,v)\\g(x)\end{pmatrix}}\,\mapsto\,\,\begin{pmatrix}
 -v\,\partial_xf\,+\,\varphi(v)\,g\\
{\displaystyle \int_{\mathbb R}\!v\,f\,dv}
\end{pmatrix},
$$
where $\,\,{\displaystyle g(x)\,=\,\frac1{4\pi e}\,E(x)}\,\,$ and
$\,\,\displaystyle{\varphi(v)\,=\,\frac{4\pi
e^2}{m}\,f_0'(v)}.\,\,$ Within the framework of stability problem
setting it is crucial to study the spectral properties of operator
$\,T\,$ acting in the underlying space prescribed by the physical
arguments
$$
\big\{{\rm L}_1(\mathbb R_x,dx)\otimes {\rm L}_1(\mathbb
R_v,(1+|\,v|)\,dv)\big\}\,\oplus\,{\rm L}_1(\mathbb R_x,dx)\,.
$$

General block operator matrix acting in the direct sum of Banach
spaces $\,X\oplus Y\,$ admits (see [2]) the following
factorization
\begin{equation}
\begin{pmatrix} A&B\\
C&D\end{pmatrix}\,\,=\,\,\begin{pmatrix} I&0\\
C\,A^{-1}&J\end{pmatrix}\begin{pmatrix} A&0\\
0&D\,-\,C\,A^{-1}B\end{pmatrix}\begin{pmatrix} I&A^{-1}B\\
0&J\end{pmatrix} \label{form7}
\end{equation}
where $\,A:X\to X\,$ is assumed to be invertible, $\,B:Y\to X\,$
and $\,C:X\to Y\,$ are bounded operators, while $\,I\,$ and
$\,J\,$ stand for identity operators in $\,X\,$ and $\,Y\,$
respectively.

\bigskip

{\bf Definition} \, {\it Expression $\,\,D\,-\,C\,A^{-1}B\,\,$ is
called the Schur complement of the block $\,A\,$ of operator
matrix {\rm (\ref{form7})}. }

\bigskip

Due to invertibility of upper and lower triangle factors in the
above formula it readily implies

\bigskip

{\bf Statement 1} \, {\it Point $\,\lambda\not\in\sigma(A)\,$
belongs to the resolvent set of operator matrix {\rm
(\ref{form7})} acting in $\,X\oplus Y\,$ if and only if  \, Schur
complement
$$
S(\lambda)\,\,:=\,\,D\,-\,\lambda\,J\,-\,C(\,A-\lambda\, I)^{-1}B
$$
is boundedly invertible in $\,Y.$
}

\bigskip

For the problem under consideration we let
$$
Y\,=\,\,{\rm L}_1(\mathbb R_x,dx)\,,\,\,\,\,X\,=\,\,Y\otimes\,
{\rm L}_1(\mathbb R_v,(1+|\,v|)\,dv)\,.
$$
In order to put our setting into the general operator theoretic
context we define the entries of block operator matrix
$$
T\,\,=\,\,\begin{pmatrix} A&B\\
C&0\end{pmatrix}
$$
as follows: \, operator
$\,A\,=\,-\,v\,\displaystyle{\frac{\partial}{\partial x}}\,\,$ is
acting in $\,X\,$ with the natural dense domain $\,\big\{f\in X\!:
f(\,\cdot,v)$ {\it absolutely continuous}, $\,v\,\partial
f/\partial x\in X\big\},\,\,$ $C\,=\,\displaystyle{\int_{\mathbb
R}\cdot\,\,v\,dv}\,\,$ and $\,\,B\,\,$ is an operator of
multiplication by function $\,\varphi(v)\,$ satisfying condition
\begin{equation}
\int_{\mathbb
R}\big(1\,+\,|v|\big)\,|\varphi(v)|\,dv\,<\,\infty\,.
\label{form8}
\end{equation}

\bigskip

{\bf Theorem 1} \, {\it The spectrum of operator $\,T\,$ consists
of two components
$$
\sigma(T)\,\,=\,\,\sigma(A)\,\cup\,\big\{\lambda\in\mathbb
C\setminus i\mathbb R: -1\in\sigma(Q(\lambda))\big\}
$$
where $\,\,\sigma(A)\,=\,i\mathbb R,\,\,$ while operator
$$
\,\,\displaystyle{Q(\lambda)\,\,=\,\,\lambda^{-1}C\,(A-\lambda\,
I)^{-1}B}
$$
acting in $\,Y\,$ is in fact a superposition of multiplication
operator $\,B:Y\to X\,$ and an integral operator
$\,\lambda^{-1}C\,(A-\lambda\, I)^{-1}\!:\,X\to Y.$ }

\bigskip

{\bf Proof} \, For $\,\lambda=\sigma+i\tau,\,\,\sigma\ne 0,\,$ the
resolvent $\,R_0(\lambda):=(A-\lambda\, I)^{-1}\,$ of operator
$\,A\,$ proves to be (see e.g. [3]) an integral operator of the
form
$$
R_0(\lambda)f(x,v)\,\,=\,\,\left\{\begin{array}{rcl}
{\displaystyle
\frac1{v}\int_x^{\infty}\exp\bigg(\frac{\lambda}{v}\,(\xi-x)\bigg)f(\xi,v)\,d\xi\,,}&&\sigma
v<0\,,\\
{\displaystyle
-\frac1{v}\int_{-\infty}^x\exp\bigg(\frac{\lambda}{v}\,(\xi-x)\bigg)f(\xi,v)\,d\xi\,,}&&\sigma
v>0\,.\\
\end{array}\right.
$$
To be specific one can assume that $\,\sigma={\rm
Re}\,\lambda>0\,$ and evaluate the norm
\begin{multline*}
\|\,R_0(\lambda)f\|_X\,\,=\,\int_{-\infty}^0(1-v)\,dv\int_{-\infty}^{\infty}dx\,\left|\,\frac1{v}
\int_{x}^{\infty}\exp\bigg(\frac{\lambda}{v}\,(\xi-x)\bigg)f(\xi,v)\,d\xi\,\right|\,\,+\\
+\,\int_0^{\infty}(1+v)\,dv\int_{-\infty}^{\infty}dx\,\left|\,\frac1{v}
\int_{-\infty}^x\exp\bigg(\frac{\lambda}{v}\,(\xi-x)\bigg)f(\xi,v)\,d\xi\,\right|\,\,\leqslant\\
\leqslant\,\int_{-\infty}^0(v-1)\,dv\int_{-\infty}^{\infty}dx\,\left(\,\frac1{v}
\int_{x}^{\infty}\exp\bigg(\frac{\sigma}{v}\,(\xi-x)\bigg)|f(\xi,v)|\,d\xi\,\right)\,\,+\\
+\,\int_0^{\infty}(1+v)\,dv\int_{-\infty}^{\infty}dx\,\left(\,\frac1{v}
\int_{-\infty}^x\exp\bigg(\frac{\sigma}{v}\,(\xi-x)\bigg)|f(\xi,v)|\,d\xi\,\right)\,\,=\\
=\,\,\frac1{\sigma}\int_{-\infty}^{\infty}(1+|\,v|)\,dv\int_{-\infty}^{\infty}|f(\xi,v)|\,d\xi\,\,=\,\,
\frac1{\sigma}\,\|f\|_X\,.
\end{multline*}
Thus we have obtained that $\,\mathbb C\setminus i\mathbb R\,$
belongs to the resolvent set of operator $\,A\,$ and moreover for
arbitrary $\lambda\not\in i\mathbb R\,$ the estimate
$$
\,\,\|\,R_0(\lambda)\,\|_{X\to X}\leqslant\,\,\,|\,{\rm
Re}\,\lambda\,|^{-1}
$$
is valid. Let us show now that spectrum of operator $\,A\,$ is
purely continuous and occupies the axis $\,\sigma(A)=i\mathbb
R.\,$ To this end given $\,\mu\in i\mathbb R\,$ it suffices to
produce in $\,X\,$ a noncompact family of approximate eigenvectors
$$
f_{\delta}(x,v)\,\,=\,\,\frac1{\delta}\,\exp\left(-\mu\,\frac{x}{v}-x^2\right)\!
\left\{\,\theta\left(\frac{v}{\delta}+1\right)-\theta\left(\frac{v}{\delta}-1\right)\right\},\,\,\,\delta>0,
$$
where $\,\theta\,$ is Heaviside step function, so that as
$\,\delta\downarrow 0\,$ one has
$$
\|f_{\delta}\|_X\,\,=\,\,\,\frac1{\delta}\int_{-\infty}^{\infty}\!dx
\int_{-\delta}^{\delta}\left|\,\exp\left(-\mu\,\frac{x}{v}-x^2\right)\right|\big(1+|v|\big)\,dv\,\,=\,\,(2+\delta)\sqrt{\pi}\,,
$$
$$
\|(A-\mu
I)f_{\delta}\|_X\,=\,\frac2{\delta}\int_{-\infty}^{\infty}\!\!|x|\,dx\!
\int_{-\delta}^{\delta}\left|\,\exp\left(-\mu\,\frac{x}{v}-x^2\right)\right|\!\big(|v|+v^2\big)dv\,=\,2\delta(1+2\delta/3)\,,
$$
and moreover
$$
\|Cf_{\delta}\|_Y\,=\,\,\frac1{\delta}\int_{-\infty}^{\infty}\!dx\,
\left|\int_{-\delta}^{\delta}\exp\left(-\mu\,\frac{x}{v}-x^2\right)v\,dv\,\right|\,\,\leqslant\,\,2\delta\,.
$$
Thus vector functions $\,(f_{\delta}(x,v),0)\,$ provides us with a
non-compact in $\,X\oplus Y\,$ family of approximate eigenvectors
for operator $\,T\,$ corresponding to the point $\,\mu\in i\mathbb
R\,$ and hence $\,i\mathbb R\subset\sigma(T).\,$

To complete the proof of Theorem 1 one should just take into
account that $\sigma(A)=i\mathbb R\,$ and make usage of Statement
1 according to which $\,\lambda\in\sigma(T)\setminus i\mathbb
R\,\,$ if and only if $\,\,-1\in\sigma(Q(\lambda))\,\,$ where
$\,\,Q(\lambda)\,=\,\lambda^{-1}CR_0(\lambda)B\,\,$ is an integral
operator given in $\,Y\,$ by a superposition formula
$$
Q(\lambda)g(x)\,\,=\,\,\frac1{\lambda}\int_{-\infty}^{\infty}v\,R_0(\lambda)\varphi(v)g(x)\,dv\,.
$$

\medskip

Theorem 1 enables us to specify and effectively localize the zone
at the complex plane $\,\mathbb C\,$ disjoint with the spectrum of
$\,T\,$ in which operator function $\,J+Q(\lambda)\,$ proves to be
invertible due to an appropriate smallness of $\,\,Q(\lambda):Y\to
Y.$

\bigskip
\smallskip

{\bf Proposition 1} \, {\it Resolvent set of operator $\,T\,$
contains the domain $\,\Omega\subset\mathbb C\,$ specified by
condition
\begin{equation}
\,\,\displaystyle{|\,{\rm
Re}\,\lambda\,|\,\,>\,\,\,\frac1{|\lambda|}\int_{-\infty}^{\infty}|\,\varphi(v)|\,|\,v|\,dv}\,.\,\,
\label{form9}
\end{equation}
}

\medskip

{\bf Proof} \, By virtue of Theorem 1 it suffices to establish the
inequality
$$
\|\,Q(\lambda)\|_{Y\to
Y}\,\,\leqslant\,\,\,\frac1{|\,\lambda|\,|\,\sigma|}\int_{-\infty}^{\infty}|\,\varphi(v)|\,|\,v|\,dv\,
$$
provided that $\,\sigma={\rm Re}\,\lambda\ne 0.\,$ To be definite
let us consider the case $\,\sigma={\rm Re}\,\lambda>0\,$ when
operator $\,Q(\lambda)\,$ is given by the expression
\begin{multline*}
Q(\lambda)g(x)\,\,=\,\,\frac1{\lambda}\,\Bigg\{\int_{-\infty}^0\varphi(v)\,dv
\int_x^{\infty}\exp\bigg(\frac{\lambda}{v}\,(\xi-x)\bigg)g(\xi)\,d\xi\,\,\,-\\
-\,\int_0^{\infty}\varphi(v)\,dv
\int_{-\infty}^x\exp\bigg(\frac{\lambda}{v}\,(\xi-x)\bigg)g(\xi)\,d\xi\,\Bigg\}.
\end{multline*}
For such $\,\lambda\,$ the following estimate
\begin{multline*}
\|Q(\lambda)g\|_Y\,\,\leqslant\,\,\frac1{|\lambda|}\,\int_{-\infty}^{\infty}dx\,\Bigg\{\int_{-\infty}^0|\,\varphi(v)|\,dv
\int_x^{\infty}\exp\bigg(\frac{\sigma}{v}\,(\xi-x)\bigg)|\,g(\xi)|\,d\xi\,\,\,+\\
+\,\int_0^{\infty}|\,\varphi(v)|\,dv
\int_{-\infty}^x\exp\bigg(\frac{\sigma}{v}\,(\xi-x)\bigg)|\,g(\xi)|\,d\xi\,\Bigg\}\,\,=\\
=\,\,
\frac1{|\,\lambda|\,\,\sigma}\,\int_{-\infty}^{\infty}|\,\varphi(v)|\,|\,v|\,dv\int_{-\infty}^{\infty}|\,g(\xi)|\,d\xi\,
\end{multline*}
holds true and similarly the case $\,\sigma={\rm Re}\,\lambda<0\,$
is dealt with.

Condition (\ref{form9}) clearly implies the inequality
$\,\|\,Q(\lambda)\|_{Y\to Y}\!<1\,$ and thus guarantees
invertibility of operator $\,J+Q(\lambda)\,$ which is necessary
and sufficient for $\,\lambda\not\in i\mathbb R\,$ to belong to
the resolvent set of operator $\,T.\,$ The boundary of the domain
$\,\Omega\,$ parametrized in variables $\,\,\sigma={\rm
Re}\,\lambda\,\,$ and $\,\,\tau={\rm Im}\,\lambda\,\,$ proves to
be two-component algebraic curve
$$
\sigma^4\,+\,\sigma^2\tau^2\,-\,c^2\,=\,0,\quad
c\,=\,\int_{-\infty}^{\infty}|\,\varphi(v)|\,|\,v|\,dv\,,
$$
intersecting the real axis at points $\,\lambda=\pm\sqrt{c}.\,$
Connected component of the boundary $\,\partial\Omega\,$ located
in the right half-plane $\,\mathbb C_+\,$ approaches the imaginary
axis at infinity with asymptotics
$\,\,\displaystyle{\sigma\,\sim\,\pm\,\frac{c}{\tau}},\,\,$ while
in the vicinity of the real axis it can be approximated by a
parabolic pattern
$\,\displaystyle{\sigma-\sqrt{c}\,\,\sim\,-\,\frac{\tau^2}{4\sqrt{c}}}.\,\,$

\medskip

{\bf Remark} \, The above information about the resolvent set of
operator $\,T\,$ associated with the problem in question specifies
and supplements localization of spectrum free zone established in
[1].

\bigskip
\bigskip

\begin{center}
\bf \S\,3. Fourier representation and Landau dispersion
relationship
\end{center}

\bigskip
\medskip

Let us denote by $\,\Phi\,$ the standard Fourier transform in the
space variable
$$
\Phi:\,\,Y\ni
f\,\,\,\longmapsto\,\,\,\widehat{f}(k)\,=\,\frac1{\sqrt{2\pi}}\int_{-\infty}^{\infty}e^{-ik
x}f(x)\,dx\,.
$$
The image $\,\widehat{Y}=\Phi\,Y\,$ will be regarded as a Banach
space equipped by the induced norm
$\,\,\|\,\widehat{f}\,\|_{\widehat{Y}}=\,\|\,f\,\|_Y\,$ with
respect to which mapping $\,\Phi\,$ clearly becomes an isometric
isomorphism.

\bigskip

{\bf Proposition 2} \, {\it After passing to Fourier
representation the normalized Schur complement
$$
-\,\lambda^{-1}\Phi\,S(\lambda)\,\Phi^{-1}\,=\,\,\Phi\,\big(J+Q(\lambda)\big)\,\Phi^{-1}
$$
associated with block matrix $\,T\,$ becomes an operator of
multiplication by the function
$$
\Delta(k,\lambda)\,\,=\,\,1\,\,+\,\,\frac1{ik}\int_{-\infty}^{\infty}\frac{\varphi(v)}
{ikv+\lambda}\,\,dv
$$
defined for $\,k=0\,$ by continuity
$$
\Delta(0,\lambda)\,\,=\,\,1\,\,-\,\,\frac1{\lambda^2}\int_{-\infty}^{\infty}v\,\varphi(v)\,dv\,.
$$
}

\smallskip

{\bf Proof} \, Provided that $\,f\,$ is taken from the domain of
operator $\,A\,$ one has
$$
\Phi(A-\lambda I)f(k)\,\,=\,\,(-ikv-\lambda)\Phi f(k)
$$
and consequently the resolvent
$\,\Phi\,R_0(\lambda)\,\Phi^{-1}\!\!:\widehat{Y}\to\widehat{Y}\,\,$
is just multiplication by the function
$\,\,-\,(ikv+\lambda)^{-1}.\,\,$ Therefore operator
$\,\,J+Q(\lambda)\,=\,-\,\lambda^{-1}S(\lambda)\,\,$ written in
Fourier representation takes the form
$$
\Phi\,\big(J+Q(\lambda)\big)\,\Phi^{-1}:\,\,\widehat{f}(k)\,\,\mapsto\,\,\left\{\,1\,\,-\,\,
\frac1{\lambda}\,\int_{-\infty}^{\infty}\frac{v\,\varphi(v)}
{ikv+\lambda}\,\,dv\,\right\}\,\widehat{f}(k)\,,
$$
where
$$
\frac1{\lambda}\,\int_{-\infty}^{\infty}\frac{v\,\varphi(v)}
{ikv+\lambda}\,\,dv
\,\,=\,\,-\,\,\frac1{ik}\,\int_{-\infty}^{\infty}\frac{\varphi(v)}
{ikv+\lambda}\,\,dv\,
$$
since
$\,\displaystyle{\int_{-\infty}^{\infty}\varphi(v)\,dv\,=\,0\,.}\,\,$
Finally let us evaluate the difference
$$
\int_{-\infty}^{\infty}\frac{v\,\varphi(v)}
{ikv+\lambda}\,\,dv\,\,-\,\,\frac1{\lambda}\,\int_{-\infty}^{\infty}\!v\,\varphi(v)\,\,dv\,\,=
\,\,-\,\,\frac{ik}{\lambda}\int_{-\infty}^{\infty}\frac{v^2\,\varphi(v)}
{ikv+\lambda}\,\,dv\,\,
$$
as $\,k\to 0.\,$ To be definite consider the case when
$\,\tau={\rm Im}\,\lambda>0\,$ so that for $\,k>0\,$ small enough
the following inequalities
\begin{eqnarray*}
\bigg|\,\frac{ik}{\lambda}\int_{-\infty}^{-k^{-1/2}}\!\!\frac{v^2\,\varphi(v)}
{ikv+\lambda}\,\,dv\,\,\bigg|&\leqslant&\frac2{|\lambda|}\,\bigg(1+\frac{\tau}{|\sigma|}\bigg)\!
\int_{-\infty}^{-k^{-1/2}}\!\!\!\!\!\!|v|\,|\varphi(v)|\,dv\,,\\
\bigg|\,\frac{ik}{\lambda}\int_{k^{-1/2}}^{-k^{-1/2}}\!\!\frac{v^2\,\varphi(v)}
{ikv+\lambda}\,\,dv\,\,\bigg|&\leqslant&\,\,\frac{\sqrt{k}}{|\lambda|\,|\sigma|}\,\int_{-\infty}^{\infty}\!|v|\,|\varphi(v)|\,dv\,,\\
\bigg|\,\frac{ik}{\lambda}\int_{k^{-1/2}}^{\infty}\frac{v^2\,\varphi(v)}
{ikv+\lambda}\,\,dv\,\,\bigg|\,\,&\leqslant&\,\,\,\frac1{|\lambda|}\,\int_{k^{-1/2}}^{\infty}\,\,|v|\,|\varphi(v)|\,dv\,,\\
\end{eqnarray*}
are valid where $\,\sigma={\rm Re}\,\lambda\ne0.\,$ Due to
condition (\ref{form8}) the right-hand sides of the above
estimates vanish as $\,k\downarrow 0\,$ and therefore
$$
\lim_{k\downarrow 0}\int_{-\infty}^{\infty}\frac{v\,\varphi(v)}
{ikv+\lambda}\,\,dv\,\,=\,\,\frac1{\lambda}\,\int_{-\infty}^{\infty}\!v\,\varphi(v)\,\,dv\,\,
$$
provided that $\,\tau={\rm Im}\,\lambda>0.\,$ The remaining cases
are dealt with similarly.

\bigskip

In physical literature (see e.g. [4]) equation
\begin{equation}
\Delta(k,\lambda)\,=\,0
\label{form10}
\end{equation}
is known to be called Landau dispersion relationship. It specifies
the values of spectral parameter $\,\lambda\,$ and wave number
$\,k\,$ for which the problem (\ref{form1})-(\ref{form2})
possesses unstable modes corresponding to undamped plasma
oscillations. The roots $\,\lambda=\lambda(k)\,$ of equation
(\ref{form10}) will be regarded as {\it singular values} of the
problem in question associated to given $\,k\in\mathbb R.$

\bigskip

{\bf Lemma 1} \, {\it Suppose that function $\,\varphi(v)\,$
satisfies the following condition
$$
\int_{-\infty}^{\infty}|\,\varphi(v)|\,|\,v|^3\,dv\,\,<\,\,\infty\,.
$$
Then for arbitrary $\,\lambda=\sigma+i\tau,\,\sigma\ne 0,\,$
function $\,\phi(k,\lambda)\,=\,\Delta(k,\lambda)-1\,$ is bounded
uniformly together with its two first derivatives in $\,k\,$ and
moreover the asymptotic estimate is valid
$$
|\,\phi(k,\lambda)\,|\,\,+\,\,|\,\phi'(k,\lambda)\,|\,\,+\,\,|\,\phi''(k,\lambda)\,|\,\,=\,\,O\big(|\,k|^{-3/2}\big)\,,\,\,
\,k\to\pm\infty.
$$ }

\medskip

{\bf Proof} \, Taking advantage of the formula $\,\,\displaystyle{
\phi(k,\lambda)\,=\,-\,\frac1{\lambda}\,\int_{-\infty}^{\infty}\frac{v\,\varphi(v)}
{ikv+\lambda}\,\,dv}\,\,$ we obtain the estimate
$$
|\,\phi(k,\lambda)|\,\,\leqslant\,\,\frac1{|\,\lambda|\,|\,\sigma|}\,\int_{-\infty}^{\infty}|\,\varphi(v)|\,|\,v|\,dv\,.
$$
To evaluate $\,\phi(k,\lambda)\,$ as $\,k\to\pm\infty\,$ it makes
sense to use representation
$$
\phi(k,\lambda)\,\,=\,\,\frac1{ik}\,\left\{\int_{-\infty}^{-|k|^{-1/2}}\!\!+\,\,\int_{-|k|^{-1/2}}^{|k|^{-1/2}}+\,\,
\int_{|k|^{-1/2}}^{\infty}\right\}\,\frac{\varphi(v)}{ikv+\lambda}\,\,dv\,,
$$
where
$$
\left|\,\int_{-|k|^{-1/2}}^{|k|^{-1/2}}\frac{\varphi(v)}{ikv+\lambda
}\,\,dv\,\,\right|\,\,\leqslant\,\,\frac2{|\,\sigma|\,\sqrt{|\,k|}}\,\,\max_{|v|\leqslant
1}|\,\varphi|
$$
and
$$
\left|\,\left\{\int_{-\infty}^{-|k|^{-1/2}}\!\!+\,\,
\int_{|k|^{-1/2}}^{\infty}\right\}\,\frac{\varphi(v)}{ikv+\lambda
}\,\,dv\,\,\right|\,\,\leqslant\,\,\frac{2\sqrt{|k|}}{|\,k|-\tau^2}\,\int_{-\infty}^{\infty}|\,\varphi(v)|\,dv
$$

\medskip

\noindent for $\,|\,k|\,$ sufficiently large. Along the same lines
one can estimate the derivatives

$$
\phi'(k,\lambda)\,=\,
\frac{i}{\lambda}\int_{-\infty}^{\infty}\!\frac{v^2\,\varphi(v)}
{(ikv+\lambda)^2}\,\,dv
\,\,=\,\,-\,\frac1{ik^2}\int_{-\infty}^{\infty}\!\frac{\varphi(v)}
{ikv+\lambda}\,dv\,\,-\,\,\frac1{k}\int_{-\infty}^{\infty}\!\frac{v\,\varphi(v)}
{(ikv+\lambda)^2}\,dv\,,
$$
\begin{multline*}
\phi''(k,\lambda)\,\,=\,\,
\frac2{\lambda}\,\int_{-\infty}^{\infty}\frac{v^3\,\varphi(v)}
{(ikv+\lambda)^3}\,\,dv
\,\,=\,\,\frac2{ik^3}\,\int_{-\infty}^{\infty}\frac{\varphi(v)}
{ikv+\lambda}\,\,dv\,\,+\\
+\,\,\frac2{k^2}\,\int_{-\infty}^{\infty}\frac{v\,\varphi(v)}
{(ikv+\lambda)^2}\,\,dv\,\,+\,\,\frac{2i}{k}\int_{-\infty}^{\infty}\frac{v^2\,\varphi(v)}
{(ikv+\lambda)^3}\,\,dv\,.
\end{multline*}

\bigskip

Specification of the resolvent set for operator $\,T\,$ in terms
of the corresponding Landau dispersion relationship is given by

\bigskip

{\bf Theorem 2} \, {\it If $\,\lambda\not\in i\mathbb R\,$ proves
not to be a root of equation {\rm(\ref{form10})} for any
$\,k\in\mathbb R\,$ then operator $\,J+Q(\lambda)\,$ is boundedly
invertible in $\,Y\,$ and moreover point $\,\lambda\,$ belongs to
the resolvent set of operator $\,T.$}

\bigskip

{\bf Proof} \,\, By Proposition 2 \,\,operator
$\,\,(J+Q(\lambda))^{-1}\,$ after passing to Fourier
representation reduces to multiplication by the function
$\,1/\Delta(k,\lambda).\,$ Due to this fact and with the
assumption $\,\,\displaystyle{\min_{k\in\mathbb
R}|\,\Delta(k,\lambda)|>0}\,\,$ taken into account in virtue of
Lemma 1 \,the functions
$$
h(k,\lambda)\,\,:=\,\,\frac1{\Delta(k,\lambda)}\,-\,1\,\,=\,\,-\,\frac{\phi(k,\lambda)}{\Delta(k,\lambda)}\,,
$$
$$
h'(k,\lambda)\,\,=\,\,-\,\frac{\phi'(k,\lambda)}{\Delta(k,\lambda)^2}\,,\,\,\,\,h''(k,\lambda)\,\,=\,\,-\,
\frac{\phi''(k,\lambda)}{\Delta(k,\lambda)^2}\,+\,2\,\frac{\phi'(k,\lambda)^2}{\Delta(k,\lambda)^3}
$$

\bigskip

\noindent are absolutely integrable in variable $\,k\,$ on the
whole axis. It follows readily that
$\,\,g(x)\,=\,\Phi^{-1}h(x)\,\in\,Y.\,$ In fact
$\,h(\pm\,\infty,\lambda)=h'(\pm\,\infty,\lambda)=0\,\,$ and hence
\begin{multline*}
\sqrt{2\pi}\,g(x)\,\,=\,\,\int_{-\infty}^{\infty}h(k,\lambda)e^{ikx}\,dk\,\,=\\
=\,\,-\frac1{ix}
\int_{-\infty}^{\infty}h'(k,\lambda)e^{ikx}\,dk\,\,=\,\,-\frac1{x^2}\int_{-\infty}^{\infty}h''(k,\lambda)e^{ikx}\,dk\,
\end{multline*}
is uniformly bounded and sufficiently rapidly decreasing at
infinity.

Therefore for arbitrary $\,f\in Y\,$ one has
\begin{multline*}
\Phi\,\big(J+Q(\lambda)\big)^{-1}\Phi^{-1}\widehat{f}(k)\,=\,\big(1\,+\,h(k,\lambda)\big)\widehat{f}(k)\,=\\
=\,\big(1\,+\,\Phi g(k)\big)\Phi f(k)\,=\,\Phi\big(f\,-\,f\!\ast
g\,\big)(k)\,,
\end{multline*}
where the convolution $\,\,\displaystyle{(f\ast
g)(x)\,=\,\int_{-\infty}^{\infty}\!f(x-y)g(y)\,dy}\,\,$ is
absolutely integrable on $\,\mathbb R\,$ so that
$$
\|f\ast g\|_Y\,\,\leqslant\,\int_{-\infty}^{\infty}
dx\int_{-\infty}^{\infty}|f(x-y)|\,|g(y)|\,dy\,\,=\,\,\|f\|_Y\,\|g\|_Y\,.
$$

\medskip

\noindent Finally we obtain the following estimate
\begin{multline*}
\big\|\big(J+Q(\lambda)\big)^{-1}f\,\big\|_Y\,=\,\,
\big\|\,\Phi\big(J+Q(\lambda)\big)^{-1}\Phi^{-1}\widehat{f}\,\,\big\|_{\widehat{Y}}\,=\\
=\,\,\big\|\Phi\big(f\,-\,f\!\ast
g\,\big)\big\|_{\widehat{Y}}\,=\,\,\big\|f\,-\,f\!\ast
g\,\big\|_Y\,\,\leqslant\,\,\big(\,1+\|g\|_Y\big)\|f\|_Y,
\end{multline*}
where
\begin{multline*}
\|\,g\|_Y\,\,=\,\,\frac1{\sqrt{2\pi}}\,\bigg(\int_{-1}^1
\bigg|\int_{-\infty}^{\infty}\!h(k,\lambda)\,e^{ikx}\,dk\,\,\bigg|\,\,dx\,\,\,+\\
+\,\,\,\bigg\{\int_{-\infty}^{-1}\!+\,\int_1^{\infty}\bigg\}\,
\bigg|\int_{-\infty}^{\infty}h''(k,\lambda)\,e^{ikx}\,dk\,\,\bigg|\,\,\frac{dx}{x^2}\,\,\bigg)\,\,\leqslant\\
\leqslant\,\,\,
\sqrt{\frac2{\pi}}\int_{-\infty}^{\infty}\Big(\,|h(k,\lambda)|\,+\,|h''(k,\lambda)|\,\Big)\,dk\,.
\end{multline*}

Thus under the hypothesis of Theorem 2 \,operator
$\,J+Q(\lambda)\,$ proves to be boundedly invertible in $\,Y\,$
and moreover
$$
\big\|\,\big(J+Q(\lambda)\big)^{-1}\big\|_{Y\to
Y}\!\leqslant\,\,\,1\,+\,\,\sqrt{\frac2{\pi}}\int_{-\infty}^{\infty}\Big(\,|h(k,\lambda)|\,+\,|h''(k,\lambda)|\,\Big)\,dk\,.
$$
To complete the proof it suffices to apply Theorem 1 \,according
to which $\,\lambda\not\in i\mathbb R\,$ belongs to the resolvent
set of operator $\,T\,$ if and only if
$\,-1\not\in\sigma(Q(\lambda)).$

\bigskip
\bigskip

\begin{center}
\bf \S\,4. Formula for instability index
\end{center}

\bigskip
\medskip

In what follows we will denote by $\,\Lambda(k)\,$ the set of
singular values of the problem in question corresponding to a
given $\,k\in\mathbb R:$
$$
\Lambda(k)\,\,:=\,\,\big\{\,\lambda\in\mathbb C\setminus i\mathbb
R:\,\Delta(k,\lambda)=0\,\big\}\,.
$$

\bigskip

{\bf Lemma 2} \, {\it The set $\,\Lambda(k)\subset\mathbb C\,$ is
symmetric with respect to imaginary axis. }

\bigskip

In fact provided that $\,\lambda=\sigma+i\tau\in\Lambda(k)\,$ one
has $\,-\overline{\lambda}=-\sigma+i\tau\in\Lambda(k)\,\,$ since
\begin{eqnarray*}
{\rm Re}\,\Delta(k,-\overline{\lambda})&=&1\,\,-\,\,\frac1{k}
\int_{-\infty}^{\infty}\frac{\varphi(v)\,(kv+\tau)}{\big((kv+\tau)^2+\sigma^2\big)^{1/2}}\,\,dv\,\,=\,\,{\rm
Re}\,\Delta(k,\lambda)\,,
\\
{\rm Im}\,\Delta(k,-\overline{\lambda})&=&\,\,\frac{\sigma}{k}\,
\int_{-\infty}^{\infty}\,\frac{\varphi(v)}{\big((kv+\tau)^2+\sigma^2\big)^{1/2}}\,\,\,dv\,\,\,=\,\,\,-\,\,{\rm
Im}\,\Delta(k,\lambda)\,.
\end{eqnarray*}

\bigskip

An effective estimate for the instability index $\,N(k),\,$ i.e.
the total number of singular values $\,\lambda=\lambda(k)\,$
corresponding to fixed $\,k\in\mathbb R\,$ and such that $\,{\rm
Re}\,\lambda(k)>0,\,$ is given by

\bigskip

{\bf Theorem 3} \, {\it Assume that $\,\varphi\in{\rm L}_1(\mathbb
R),\,\,\varphi(v)\to 0\,$ as $\,v\to\pm\infty,\,$ and the
derivative $\,\varphi'(v)\,$ is bounded. Every zero $\,s\,$ of
function $\,\varphi\,$ is supposed to be non-degenerate and each
of them satisfy the following condition
\begin{equation}
1\,\,-\,\,\frac1{k^2}\int_{-\infty}^{\infty}\frac{\varphi(v)}{v-s}\,dv\,\,\ne\,\,0\,
\end{equation}
for a certain $\,k\in\mathbb R\setminus\{0\}.\,$ Then the total
number $\,N(k)\,$ of the roots of equation {\rm(\ref{form10})}
located in the right half-plane are evaluated by the formula
$$
N(k)\,\,:=\,\,\#\,\big\{\lambda\in\Lambda(k)\!:\,{\rm
Re}\,\lambda>0\big\}\,\,=\,\,N_+(k)\,-\,N_-(k)\,,
$$
where
\begin{eqnarray*}
N_+(k)&=&\#\,\bigg\{s\in\mathbb
R\!:\,\varphi(s)=0,\,\,\varphi'(s)>0,\,\,1\,-\,\frac1{k^2}\int_{-\infty}^{\infty}\frac{\varphi(v)}{v-s}\,dv\,<\,0\bigg\},\\
N_-(k)&=&\#\,\bigg\{s\in\mathbb
R\!:\,\varphi(s)=0,\,\,\varphi'(s)<0,\,\,1\,-\,\frac1{k^2}\int_{-\infty}^{\infty}\frac{\varphi(v)}{v-s}\,dv\,<\,0\bigg\}.
\end{eqnarray*}
The right half-plane contains at most one root of equation
$\,\Delta(0,\lambda)=0\,$ and, moreover, $\,\,\Lambda(0)\cap\{{\rm
Re}\,\lambda>0\}\ne\varnothing\,\,$ if and only if
$$
\int_{-\infty}^{\infty}v\,\varphi(v)\,dv\,\,>\,\,0\,.
$$
}

\medskip

{\bf Proof} \, Let us fix $\,k\in\mathbb R\setminus\{0\},\,$
introduce notation $\,\,\psi(v)\,:=\,\varphi(v)/k^2\,\,$ and
consider the function
$$
w(z)\,\,=\,\,1\,\,-\,\int_{-\infty}^{\infty}\frac{\psi(v)}{v-z}\,\,dv\,\,=\,\,\Delta(k,-ikz)\,
$$
analytic in $\,\mathbb C_{\pm}\,$ whose boundary values at the
real axis are calculated by Sokhotski-Plemelj formulas
$$
w(s\pm i0)\,\,=\,\,1\,\,-\,\,{\rm
v.p.}\int_{-\infty}^{\infty}\frac{\psi(v)}{v-s}\,\,dv\,\,\mp\,\,
i\pi\psi(s)\,.
$$

\medskip

In order to evaluate the total number of zeroes of the function
$\,w(z)\,$ in the upper half-plane $\,\mathbb C_+\,$ which
coincides with the required quantity
$\,\,\#\,\big\{\lambda\in\Lambda(k)\!:\,{\rm
Re}\,\lambda>0\big\}\,\,$ we will take advantage of the argument
principle (see e.g. [5]). Beforehand let us verify that $\,w(z)\to
1\,$ as $\,\mathbb C_+\ni z\to\infty.\,$ In fact for
$\,z=s+it\in\mathbb C_+\,$ one has
$$
|\,w(z)-1|\,\,=\,\,\left|\,\int_{-\infty}^{\infty}\frac{\psi(v)}{v-s-it}\,\,dv\,\right|\,\,\leqslant\,\,
\frac1{t}\int_{-\infty}^{\infty}|\,\psi(v)|\,dv\,\,\to\,\,0\,
$$
as $\,t={\rm Im}\,z\to\infty.\,$ To estimate absolute value
$\,|\,w(z)-1|\,$ in the case when $\,s={\rm Re}\,z\to+\infty\,$ it
makes sense to represent the corresponding integral in the form
$$
1-w(z)\,\,=\,\,\left\{\int_{-\infty}^{s/2}\!+\,\,\int_{s/2}^{s-\varepsilon}\!+\,\,
\int_{s-\varepsilon}^{s+\varepsilon}\!+\,\,\int_{s+\varepsilon}^{\infty}\right\}\frac{\psi(v)}{v-s-it}\,\,dv\,,
$$
where $\,\varepsilon\in(0,s/2),\,$ and treat the above summands
separately\,:
\begin{multline*}
\left|\,\int_{s-\varepsilon}^{s+\varepsilon}\frac{\psi(v)}{v-s-it}\,\,dv\,\,-
\,\,\psi(s)\big(\,\pi\,-\,2\,\arctan\big(t/\varepsilon\big)\big)\,\right|\,\,=\\
=\,\,\,\left|\,\int_{s-\varepsilon}^{s+\varepsilon}\frac{\psi(v)-\psi(s)}{v-s-it}\,\,dv\,\right|\,\,\leqslant\,\,
2\,\varepsilon\,\max|\,\psi'(v)|\,
\end{multline*}
and, besides,
$$
\left|\,\int_{-\infty}^{s/2}\frac{\psi(v)}{v-s-it}\,\,dv\,\right|\,\,\leqslant\,\,
\frac2{s}\int_{-\infty}^{\infty}|\,\psi(v)|\,dv\,\,\to\,\,0\,,
$$
$$
\left|\,\left\{\int_{s/2}^{s-\varepsilon}\!+\,\,\int_{s+\varepsilon}^{\infty}\right\}
\frac{\psi(v)}{v-s-it}\,\,dv\,\right|\,\,\leqslant\,\,
\frac1{\varepsilon}\int_{s/2}^{\infty}|\,\psi(v)|\,dv\,\,\to\,\,0\,.
$$
\vspace{0.1cm}

\noindent Thus the following inequality
$$
|\,w(z)-1\,|\,\,\leqslant\,\,2\,\varepsilon\,\max|\,\psi'(v)|\,\,+\,\,\pi\,|\psi(s)|
\,\,+\,\,\frac2{s}\int_{-\infty}^{\infty}|\psi(v)|\,dv\,\,+\,\,
\frac1{\varepsilon}\int_{s/2}^{\infty}|\psi(v)|\,dv\,
$$
holds true where the first summand on the right-hand side can be
made arbitrarily small under an appropriate choice of
$\,\varepsilon\in(0,s/2)\,$ while the second, the third and the
fourth ones vanish for fixed $\,\varepsilon>0\,$ as
$\,s\to+\infty.\,$ Similarly one can carry out estimation of the
absolute value $\,|\,w(z)-1|\,$ in the case when $\,s={\rm
Re}\,z\to-\infty.\,$

\medskip

According to Cauchy argument principle given arbitrary closed
contour $\,\gamma\subset\mathbb C_+\,$ encircling all the zeroes
of function $\,w(z)\,$ their total multiplicity $\,N=N(k)\,$ is
equal to the sum of logarithmic residues of $\,w(z)\,$ associated
with the interior of $\,\gamma.\,$ It coincides with the index of
the point $\,w=0\,$ with respect to the curve
$\,\Gamma=w(\gamma)\,$ also known as winding number, i.e. the
total number of times that curve $\,\Gamma\,$ travels
counterclockwise around the origin:
$$
N\,\,=\,\,{\rm Ind}_{\Gamma}(0)\,\,=\,\,\frac1{2\pi
i}\oint_{\gamma}\,\,\frac{w'(z)}{w(z)}\,\,dz\,\,=\,\,\frac1{2\pi}
\oint_{\gamma}\,\,d\,\arg\,w(z)\,.
$$
Such a method of counting the roots of dispersion relationship is
related to the stability criterion due to H.\,Nyquist (see [4])
which has been rigorously justified in the context of plasma
stability problem by O.\,Penrose in [6].

\medskip

In the present setting we choose contour $\,\gamma\,$ to be
composed of the segment $\,[-R,R\,]\,$ and the half-circle
$\,z=R\,e^{i\theta},\,\theta\in[0,\pi].\,$ Radius $\,R>0\,$ is to
be taken large enough so that according to the above calculations
values $\,w(Re^{i\theta})\,$ would belong to sufficiently small
neighborhood of the point $\,w=1\,$ being separated away from the
origin. At the same time the image $\,w([-R,R\,])\,$ proves to be
a path with parametrization
\begin{eqnarray*}
{\rm Re}\,w(s+i0)&=&1\,-\,{\rm
v.p.}\int_{-\infty}^{\infty}\frac{\psi(v)}{v-s}\,dv\,,\\
{\rm Im}\,w(s+i0)&=&-\,\pi\psi(s),\,\,\,s\in[-R,R\,]\,,
\end{eqnarray*}
and such that its endpoints corresponding to $\,s=\pm R\,$ belong
to the prescribed neighborhood of $\,w=1.\,$ Moreover curve
$\,\Gamma=w(\gamma)\,$ intersects the negative semiaxis $\,\mathbb
R_-\,$ at points specified by the values of parameter
$\,s_j\in[-R,R\,]\,$ enumerated in the ascending order and such
that
$$
\psi(s_j)\,=\,0\,,\quad
1\,\,-\,\int_{-\infty}^{\infty}\frac{\psi(v)}{v-s_j}\,dv\,\,<\,\,0\,.
$$

\medskip

The gradient of transversal intersection corresponding to the zero
$\,s=s_j\,$ of function $\,\psi\,$ is determined by the derivative
$\,\psi'(s_j)\ne 0\,$ so that the arc $\,w([s_j,s_{j+1}])\,$
between two subsequent intersections of $\,\Gamma\,$ with
$\,\mathbb R_-\,$ due to condition (\ref{form9}) produce the
following increment of the argument
$$
\int_{s_j}^{s_{j+1}}\!d\,\arg\,w(s+i0)\,\,=\,\,\,\pi\,\big(\,{\rm
sign}\,\psi'(s_j)\,+\,{\rm sign}\,\psi'(s_{j+1})\big)\,.
$$
In this way the winding number for the curve $\,\Gamma\,$ is given
by the expression
\begin{equation*}
N(k)\,\,=\,\,{\rm
Ind}_{\Gamma}(0)\,\,=\,\,\int_{-\infty}^{\infty}\!d\,\arg\,w(s+i0)\,\,=
\end{equation*}
\begin{multline*}
=\,\,\frac12\,\,{\rm
sign}\,\psi'(s_1)\,\,+\,\,\,\frac12\sum_{j=1}^{n-1} \big(\,{\rm
sign}\,\psi'(s_j)\,+\,{\rm sign}\,\psi'(s_{j+1})\big)\,\,\,+\\
+\,\,\frac12\,\,{\rm
sign}\,\psi'(s_n)\,\,=\,\,\,N_+(k)\,-\,N_-(k)\,.
\end{multline*}

\bigskip

{\bf Remark} \, Singular values $\,\lambda=\lambda(k)\,$
corresponding to a fixed $\,k\ne 0\,$ can be located at the
imaginary axis $\,i\mathbb R\,$ if and only if
$$
\psi\big(i\lambda/k\big)\,=\,0\,,\quad
1\,\,-\,\int_{-\infty}^{\infty}\frac{\psi(v)}{v-i\lambda/k}\,dv\,\,=\,\,0\,.
$$
Thus condition (\ref{form9}) means that operator $\,T\,$ does not
have any singular values embedded into its continuous spectrum.

\bigskip
\medskip

{\bf Lemma 4}  \, {\it Under hypotheses of Theorem {\rm 3}
\,condition
\begin{equation}
|\,k|^4\,\,>\,\,8\,\max|\,\varphi'|\,\int_{-\infty}^{\infty}|\,\varphi(v)|\,dv
\end{equation}
guarantees absence of singular values $\,\lambda=\lambda(k)\,$ in
the right half-plane\,{\rm:}
$$\,\Lambda(k)\cap\big\{{\rm
Re}\,\lambda>0\big\}\,=\,\varnothing.
$$ }

\medskip

In fact provided that $\,\varphi(s)=0\,$ for arbitrary
$\,\varepsilon>0\,$ one has
\begin{multline*}
\left|\,\int_{-\infty}^{\infty}\frac{\varphi(v)}{v-s}\,dv\,\right|\,\,\leqslant\,\,
\left|\,\int_{s-\varepsilon}^{s+\varepsilon}\frac{\varphi(v)-\varphi(s)}{v-s}\,dv\,\right|\,\,+
\,\,\left|\,\left\{\int_{-\infty}^{s-\varepsilon}\!+\,\int_{s+\varepsilon}^{\infty}\right\}\frac{\varphi(v)}{v-s}\,dv\,\right|\,\,\leqslant\\
\leqslant\,\,
2\varepsilon\,\max|\,\varphi'(v)|\,\,+\,\,\frac1{\varepsilon}\int_{-\infty}^{\infty}|\,\varphi(v)|\,dv\,.
\end{multline*}
An optimal choice of parameter $\,\varepsilon>0\,$ minimizing the
right-hand side of the above estimate implies
$$
\left|\,\int_{-\infty}^{\infty}\frac{\varphi(v)}{v-s}\,dv\,\right|\,\,\leqslant\,\,
2\sqrt{2}\,\big(\max|\,\varphi'(v)|\,\big)^{1/2}\!\left(\int_{-\infty}^{\infty}|\,\varphi(v)|\,dv\right)^{1/2}\!.
$$
Hence for arbitrary $\,k\in\mathbb R\,$ satisfying condition
(\ref{form9}) and any zero $\,v=s\,$ of function $\,\varphi\,$ the
following inequality
$$
1\,\,-\,\,\frac1{k^2}\int_{-\infty}^{\infty}\frac{\varphi(v)}{v-s}\,dv\,\,>\,\,0\,
$$
is valid and therefore $\,N_+(k)=N_-(k)=0\,$ so that $\,
\,\#\,\big\{\lambda\in\Lambda(k)\!:\,{\rm
Re}\,\lambda>0\big\}\,=\,0\,$ by virtue of Theorem 3.

\bigskip
\bigskip

\begin{center}
\bf \S\,5. The model of two-stream instability
\end{center}

\bigskip
\medskip

A simple sufficient condition is known (see [1],[4]) to be
formulated in terms of unperturbed distribution function
$\,f_0(v)\,$ which forbids existence of unstable plasma
oscillatory perturbations. Namely provided that $\,f_0(v)\,$
possesses a unique extremum (maximum) $\,v=a\,$ the problem under
consideration has no singular values in the right half-plane:
$\,\Lambda(k)\cap\{{\rm Re}\,\lambda>0\}=\varnothing\,$ for
arbitrary $\,k\in\mathbb R.\,$ In fact one has
$\,\psi(v)/(v-a)<0\,$ and hence $\,N_-(k)=N_+(k)=0\,$ for $\,k\ne
0\,$ since
$$
1\,\,-\,\,\int_{-\infty}^{\infty}\frac{\psi(v)}{v-a}\,dv\,\,>\,\,0\,.
$$
Moreover
$\,\displaystyle{\int_{-\infty}^{\infty}\!(v-a)\varphi(v)\,dv\leqslant
0}\,\,\,$ so that $\,\Lambda(0)\cap\{{\rm
Re}\,\lambda>0\}=\varnothing\,$ by Theorem 3.

In the situation when $\,f_0(v)\,$ has two maxima a phenomenon may
happen called two-stream instability. To this end critical points
$\,v=a\,$ and $\,v=b\,$ corresponding to maximal values of the
unperturbed distribution function are to be situated so that
condition
$$
\int_a^b\frac{\psi(v)}{v-c}\,dv\,\,>\,\,1\,\,+\,\,\bigg
\{\int_{-\infty}^a\!+\,
\int_b^{\infty}\bigg\}\,\frac{\psi(v)}{c-v}\,dv
$$
holds true where $\,v=c\in(a,b)\,$ is the critical point of
$\,f_0(v)\,$ corresponding to its minimum and moreover the higher
and the wider maximum peaks should be the farther they are located
(cf. Lemmas 5 and 6 below).

To study the effect of two-stream instability we will consider the
case when all the hypotheses of Theorem 3 are satisfied and
moreover $\,f_0(a)=f_0(b)=M\,$ while $\,f_0(c)=0\,$ so that
$$
\int_{-\infty}^a \psi(v)\,dv\,\,=\,\int_c^a
\psi(v)\,dv\,\,=\,\int_c^b \psi(v)\,dv\,\,=\,\int_{\infty}^b
\psi(v)\,dv\,
$$
where $\,\psi(v)>0\,$ for $\,v\in(-\infty,a)\cup(c,b)\,$ and
$\,\psi(v)<0\,$ when $\,v\in(a,c)\cup(b,\infty)\,$ respectively.
As a consequence one has
$$
\int_{-\infty}^{\infty}
\frac{\psi(v)}{v-a}\,\,dv\,\,=\,\int_{-\infty}^a
\frac{\psi(v)}{v-a}\,\,dv\,+\,\,\bigg\{\!\int_a^c+\,\int_c^b\bigg\}
\frac{\psi(v)}{v-a}\,\,dv\,+\,\int_b^{\infty}
\frac{\psi(v)}{v-a}\,\,dv\,\,<\,\,0
$$
and similarly
$$
\int_{-\infty}^{\infty}
\frac{\psi(v)}{v-b}\,\,dv\,=\,\int_{-\infty}^a
\frac{\psi(v)}{v-b}\,\,dv\,+\,\bigg\{\!\int_a^c+\int_c^b\bigg\}
\frac{\psi(v)}{v-b}\,\,dv\,+\,\int_b^{\infty}
\frac{\psi(v)}{v-b}\,\,dv\,\,<\,\,0\,
$$
so that $\,N_-(k)=0\,$ and $\,N_+(k)\leqslant 1.\,$ A criterion of
two-stream instability within the present setting is as follows

\bigskip

{\bf Statement 2} \, {\it Given $\,k\in\mathbb R\setminus\{0\}\,$
one has $\,N_+(k)=1\,$ or equivalently
$\,\,\Lambda(k)\cap\big\{{\rm
Re}\,\lambda>0\big\}\,\ne\,\varnothing\,\,$ if and only if
$$
\int_{-\infty}^{\infty}
\frac{\varphi(v)}{v-c}\,\,dv\,\,>\,\,k^2\,.
$$
}

\bigskip

In order to formulate sufficient conditions which guarantee
existence of unstable eigenmodes in different terms we introduce
additional notation. For the critical point $\,v=a\,$ given
$\,\mu\in(0,M)\,$ set
\begin{eqnarray*}
a_<(\mu)&=&\sup\big\{\,v<a:\,f_0(v)\,\leqslant\,\mu\,\big\}\,,\\
a_>(\mu)&=&\,\inf\,\big\{\,v>a:\,f_0(v)\,\leqslant\,\mu\,\big\}\,,
\end{eqnarray*}
so that quantity $\,\big(\,a_>(\mu)\,-\,a_<(\mu)\,\big)\,$ is the
width of corresponding maximum peak at level $\,\mu,\,$ i.e. a
diameter of the connected component of preimage
$\,f_0^{-1}\big([\,\mu,\infty)\big)\,$ containing the point
$\,v=a.\,$ Similarly the values $\,b_<(\mu)\,$ and $\,b_>(\mu)\,$
are defined for the critical point $\,v=b.\,$

\bigskip

{\bf Lemma 5} \, {\it Let the inequality
$$
\frac{\big(\,a_>-a_<\big)\,\xi}{\big(\,c-a_>\big)\big(\,c-a_<\big)}\,\,+\,\,
\frac{\big(\,b_>-b_<\big)\,\eta}{\big(\,b_<-c\,\big)\big(\,b_>-c\,\big)}\,\,>\,\,
k^2\,
$$
be valid for certain $\,\xi,\eta\in(0,M),\,$ where
$\,a_>\!=a_>(\xi),\,\,a_<\!=a_<(\xi)\,$ and
$\,b_>\!=b_>(\eta),\,\,b_<\!=b_<(\eta).\,$ Then $\,\,N_+(k)=1\,\,$
so that
$$
\#\,\big\{\lambda\in\Lambda(k)\!:\,{\rm
Re}\,\lambda>0\big\}\,=\,1.\,
$$}

\medskip

In fact one has
$\,\displaystyle{\int_{a_<}^a\!\psi(v)\,dv\,=\int_{a_>}^a\!\psi(v)\,dv},\,$
hence
$\,\displaystyle{\int_{a_<}^{a_>}\frac{\psi(v)}{v-c}\,dv\,>\,0}\,\,$
and therefore
\begin{multline*}
\int_{-\infty}^c\frac{\psi(v)}{v-c}\,dv\,\,\,>\,\,\,
\left\{\int_{-\infty}^{a_<}+\,\,\int_{a_>}^c\right\}\,\frac{\psi(v)}{v-c}\,dv\,\,\,>\\
>\,\,\frac1{a_<-c}\,\int_{-\infty}^{a_<}\psi(v)\,dv\,\,+\,\,\frac1{a_>-c}\,\int_{a_>}^c\psi(v)\,dv\,\,
=\,\,\frac{\big(\,a_>-a_<\big)\,\xi/k^2}{\big(\,c-a_>\big)\big(\,c-a_<\big)}\,
\end{multline*}
since
$\,\,\displaystyle{\int_{-\infty}^{a_<}\!\!\psi(v)\,dv\,=\,\int_c^{a_>}\!\!\psi(v)\,dv\,=\,\xi/k^2.}\,\,$
Similarly the inequality
\begin{multline*}
\int_c^{\infty}\frac{\psi(v)}{v-c}\,dv\,\,\,>\,\,\,
\left\{\int_c^{b_<}+\,\,\int_{b_>}^{\infty}\right\}\,\frac{\psi(v)}{v-c}\,dv\,\,\,>\\
>\,\,\frac1{b_<-c}\,\int_c^{b_<}\psi(v)\,dv\,\,+\,\,\frac1{b_>-c}\,\int_{b_>}^{\infty}\psi(v)\,dv\,\,
=\,\,\frac{\big(\,b_>-b_<\big)\,\eta/k^2}{\big(\,b_<-c\,\big)\big(\,b_>-c\,\big)}\,
\end{multline*}
is verified where
$\,\,\displaystyle{\int_c^{b_<}\!\!\psi(v)\,dv\,=\,\int_{\infty}^{b_>}\!\!\psi(v)\,dv\,=\,\eta/k^2.}\,\,$
To complete the proof of Lemma it suffices to apply Statement 2. A
somewhat different type condition of instability is given by
\bigskip

{\bf Lemma 6} \, {\it If there exist $\,\sigma\in(0,c-a)\,$ and
$\,\tau\in(0,b-c)\,$ such that the inequality
$$
\frac{\sigma\,f_0(a+\sigma)}{\big(c-a-\sigma\big)\big(c-a\big)}\,\,+\,\,
\frac{\tau\,f_0(b-\tau)}{\big(b-c-\tau\big)\big(b-c\big)}\,\,>\,\,
k^2\,
$$
holds true then the right half-plane $\,{\rm Re}\,\lambda>0\,$
contains just one singular value $\,\lambda=\lambda(k)\,$
corresponding to given $\,k\in\mathbb R\setminus\{0\}.$}

\bigskip

Really taking into account that
$\,\,\displaystyle{\int_{-\infty}^a\psi(v)\,dv\,=\,\int_c^a\psi(v)\,dv}\,\,\,$
one has
$$
\int_{-\infty}^c\frac{\psi(v)}{v-c}\,dv\,\,>\,\,\frac1{a-c}\,\int_{-\infty}^a\!\!\psi(v)\,dv\,\,\,+\,\,
\int_a^c\frac{\psi(v)}{v-c}\,dv\,\,=\,\,\int_a^c\bigg(\frac1{v-c}\,-\,\frac1{a-c}\bigg)\psi(v)\,dv\,
$$
where $\,\psi(v)<0\,$ for $\,v\in(a,c)\,$ and hence
\begin{multline*}
\int_{-\infty}^c\frac{\psi(v)}{v-c}\,dv\,\,>\,\,\int_{a+\sigma}^c\bigg(\frac1{v-c}\,-\,\frac1{a-c}\bigg)\psi(v)\,dv\,\,\geqslant\\
\geqslant\,\,\bigg(\frac1{c-a-\sigma}\,-\,\frac1{c-a}\bigg)\int_c^{a+\sigma}\!\!\psi(v)\,dv\,\,=\,\,
\frac{\sigma\,f_0(a+\sigma)/k^2}{(c-a-\sigma)(c-a)}\,.
\end{multline*}
Along the same lines we obtain the inequality
\begin{multline*}
\int_c^{\infty}\frac{\psi(v)}{v-c}\,dv\,\,>\,\,\int_c^{b-\tau}\bigg(\frac1{v-c}\,-\,\frac1{b-c}\bigg)\psi(v)\,dv\,\,\geqslant\\
\geqslant\,\,\bigg(\frac1{b-c-\tau}\,-\,\frac1{b-c}\bigg)\int_c^{b-\tau}\!\!\psi(v)\,dv\,\,=\,\,
\frac{\tau\,f_0(b-\tau)/k^2}{(b-c-\tau)(b-c)}\,
\end{multline*}
and thus the required sufficient condition proves to be a
straightforward corollary of Statement 2.

\bigskip
\bigskip
\bigskip

\begin{center}
\bf REFERENCES
\end{center}

\bigskip
\medskip

\begin{enumerate}
\item V. P. Maslov, M. V. Fedoryuk \, The linear theory of Landau
damping, Sbornik Mathematics, 1986, v. 55, № 2, p. 437-465. \item
P. R. Halmos \, A Hilbert space problem book, Springer-Verlag,
1982. \item T. Kato \, Perturbation theory for linear operators \,
Springer-Verlag, 1966. \item T. H. Stix \, Waves in plasmas, New
York, American Institute of Physics, 1992. \item M. A. Lavrentiev,
B. V. Shabat \, Methods of the theory of functions of the complex
variable, Moscow, Nauka, 1973. \item Penrose O. \, Electrostatic
instabilities of a uniform non-Maxwellian plasma, Physics of
Fluids, 1960, v. 2, № 2, p. 258-264.
\end{enumerate}

\end{document}